# Multiresolution ORKA: fast and resolution independent object reconstruction using a K-approximation graph

Florian Bossmann[1] and Wenze Wu[1]

[1] Harbin Institute of Technology, School of Mathematics, 150001 Harbin, China
`f.bossmann@hit.edu.cn`

**Abstract.** Object recognition and reconstruction is of great interest in many research fields. Detecting pedestrians or cars in traffic cameras or tracking seismic waves in geophysical exploration are only two of many applications. Recently, the authors developed a new method – Object reconstruction using K-approximation (ORKA) – to extract such objects out of given data. In this method a special object model is used where the movement and deformation of the object can be controlled to fit the application.
ORKA in its current form is highly dependent on the data resolution. On the one hand, the movement of the object can only be reconstructed on a grid that depends on the data resolution. On the other hand, the runtime increases exponentially with the resolution. Hence, the resolution of the data needs to be in a small range where the reconstruction is accurate enough but the runtime is not too high. In this work, we present a multiresolution approach, where we combine ORKA with a wavelet decomposition of the data. The object is then reconstructed iteratively what drastically reduces the runtime. Moreover, we can increase the data resolution such that the movement reconstruction no longer depends on the original grid. We also give a brief introduction on the original ORKA algorithm. Hence, knowledge of the previous work is not required.

## 1   Introduction

Object reconstruction tries to detect or extract objects from given data. This problem occurs in many applications, such as street sign recognition for autonomous vehicles [1], subsurface reservoir detection in geophysical exploration [2], and text recognition for hand written documents [3]. The authors proposed a new algorithm for object reconstruction in an earlier work [4]. The ORKA algorithm – Object reconstruction using K-approximation - is a greedy method that uses a specially designed object model. The model allows the object to change its position and form between two measurements. The amount of change can be limited based on the physical limitations of the object. ORKA assumes that the data only contains few objects, i.e., the data is sparse in the object domain. The algorithm recovers one object at a time and repeats this process until all objects are found. However, the object model of ORKA is more complicated

2than models of other techniques. Hence, it requires a more complex algorithm to reconstruct an object. Indeed, the method proposed in [4] has a complexity of $O((2C+1)^K)$ where $K$ is a parameter controlling the approximation rate and $C$ depends on the data resolution. Moreover, the object movement is restricted to points on a grid which also depends on the data resolution. This causes a dilemma: If we want good reconstruction results, we need to choose a large parameter $K$ and a high resolution for the data, i.e., a large parameter $C$. In this work, we present a method to remove the dependence of $C$ on the resolution. The approach can further be used to refine the movement grid making it independent from the original data resolution. This finally allows to run ORKA on any data without having to balance between runtime and data resolution.

The remainder of the work is organized as follows. In the next section we introduce the notation and operators used throughout this work. We also define the object model used by ORKA and give a detailed summary on the previous work [4] that contains all results that are necessary for this work. In Section 3 we discuss the new multiresolution approach that overcomes the problems mentioned above. We also give theoretical error bounds for the new method. Finally, in the last section we give numerical examples and compare the new approach to the standard ORKA method.

## 2    Preliminaries

Throughout this work the entries of a matrix or vector are addressed using lower indexing, e.g., $\boldsymbol{A}_{jk}$ or $\boldsymbol{u}_j$. If we refer to a complete column of a matrix, we replace its index with a colon, for example $\boldsymbol{A}_{:k}$ is the $k$-th column of the matrix $\boldsymbol{A}$. We are using the Euclidean norm $\|\boldsymbol{u}\|_2$ for vectors and the Frobenius norm $\|\boldsymbol{A}\|_F$ for matrices. The corresponding inner product is denoted by $\langle \boldsymbol{A}, \boldsymbol{B} \rangle$ and $\langle \boldsymbol{u}, \boldsymbol{v} \rangle$.

The given input data is always written as $\boldsymbol{D}$. From now on, assume that $\boldsymbol{D} \in \mathbb{R}^{M \times N}$ is a matrix containing data from $N$ different measurements, where each measurement is one-dimensional. This holds for e.g., seismic data where each column in the data represents a time recording of one sensor (Fig. 2). The presented algorithm can also be applied on tensor data $\boldsymbol{D} \in \mathbb{R}^{M_1 \times M_2 \times N}$ where each measurement has two (or more) dimensions, such as video data. However, the higher dimensional case does not give any more insight to the problem. Hence, we stick to the matrix case in this work to keep the notation simple.

### 2.1    Object model and ORKA

To introduce our object model, we first need to define the column shift operator. Therefore, we need the shift matrix

$$\boldsymbol{J} = \begin{pmatrix} 0 & \cdots & 0 & 1 \\ 1 & \ddots & \vdots & 0 \\ 0 & \ddots & \ddots & \vdots \\ \ddots & 0 & 1 & 0 \end{pmatrix}.$$

**Definition 1.** For an integer vector $\boldsymbol{\lambda} \in \mathbb{N}^N$, define $S_{\boldsymbol{\lambda}} \colon \mathbb{R}^{M \times N} \to \mathbb{R}^{M \times N}$ as

$$S_\lambda(A) = \left(J^{\lambda_1} A_{:1}, J^{\lambda_2} A_{:2}, \cdots, J^{\lambda_N} A_{:N}\right).$$

The operator shifts the $k$-th column of a matrix by $\lambda_k$ entries. For vectors $\boldsymbol{\lambda}, \boldsymbol{\mu} \in \mathbb{N}^N$ and matrices $\boldsymbol{A}, \boldsymbol{B} \in \mathbb{R}^{M \times N}$ the following properties follow directly from the definition:

$$\begin{aligned}
S_\lambda(A+B) = S_\lambda(A) + S_\lambda(B), \quad & S_\lambda\left(S_\mu(A)\right) = S_{\lambda+\mu}(A), \quad S_\lambda^{-1} = S_{-\lambda}, \\
\langle S_\lambda(A), S_\lambda(B)\rangle = \langle A, B\rangle, \quad & \|S_\lambda(A)\|_F = \|A\|_F.
\end{aligned} \quad (1)$$

The shift operator is used to represent the object movement in our model. Each column of $\boldsymbol{D}$ corresponds to a different observation or measurement. Hence, with the above column shift operator we can move an object into different positions for each measurement. The complete model is given by the following definition.

**Definition 2.** Let $\boldsymbol{U} \in \mathbb{R}^{M \times N}$ and $\boldsymbol{\lambda} \in \mathbb{N}^N$ be given. We call the matrix $S_\lambda(\boldsymbol{U})$ an object. We say that it moves Lipschitz-continuous with constant $C$, if $|\lambda_k - \lambda_{k+1}| \leq C$ holds for all $k = 1, \ldots, N-1$. Furthermore, we measure its total change by

$$\sum_{k=1}^{N-1} \left\|\boldsymbol{U}_{:k} - \boldsymbol{U}_{:(k+1)}\right\|_2^2. \tag{2}$$

The above definition is a simplification for equidistantly sampled data, i.e., we assume that the $N$ measurements were taken from positions $x_k$ on an equidistant grid. For non-equidistant samplings, the Lipschitz inequality needs to consider the different distances $|\lambda_k - \lambda_{k+1}| \leq C |x_k - x_{k+1}|$, and in the total change each summand must be weighted by the inverse distance $|x_k - x_{k+1}|^{-1}$. We stick to equidistant sampling in this work to simplify the notation.

The model is motivated as follows. Each column $\boldsymbol{U}_{:k}$ represents the objects appearance as it can be observed in the $k$-th measurement, while the shift $\lambda_k$ controls its position. Both, position and appearance, are subject to physical limitations and do not change arbitrarily from one measurement to the other. For the movement of the object, we assume that it e.g., has a maximum speed that can never be broken. This directly results in the shift vector $\boldsymbol{\lambda}$ being Lipschitz-continuous. The change in appearance, e.g., due to deformation or rotation, can be more complicated. Consider a car in a video that drives away from the camera and then takes a left turn. Within a few frames the car view will change from its rear to a side view. Hence, the change in appearance of an object can be strong from one frame to the other. However, we assume that these effects are rare and the total change as measured with Eq. (2) still stays small.

Now we have everything at hand to introduce the ORKA algorithm. The idea of the algorithm is, to find the object that best matches the given data, i.e., to solve

$$\min_{\lambda, U} \|D - S_\lambda(U)\|_F^2 + \mu \sum_{k=1}^{N-1} \left\|\boldsymbol{U}_{:k} - \boldsymbol{U}_{:(k+1)}\right\|_2^2 \quad \text{s.t.} \quad |\lambda_k - \lambda_{k+1}| \leq C. \tag{3}$$

The first term is the data fidelity term and fits the object to our measurements. The second term penalizes the total change of an object. The parameter $\mu$ balances both terms and should be chosen according to the noise level of the data. The continuity of the movement is added as strict constraints as it should never be violated.

Using Eq. (1) we can rewrite the data fidelity term in Eq. (3) as

$$\|\boldsymbol{D} - S_\lambda(\boldsymbol{U})\|_F^2 = \|S_{-\lambda}(\boldsymbol{D} - S_\lambda(\boldsymbol{U}))\|_F^2 = \|S_{-\lambda}(\boldsymbol{D}) - \boldsymbol{U}\|_F^2.$$

This separates the variables $\boldsymbol{\lambda}$ and $\boldsymbol{U}$. For fixed $\boldsymbol{\lambda}$, Eq. (3) is a quadratic minimization problem in $\boldsymbol{U}$. Moreover, it is convex as all involved functions are convex norms. Hence, we can derive its minimum analytically. The minimal value is

$$-\langle \boldsymbol{A}^{-1}, (S_{-\lambda}(\boldsymbol{D}))^T S_{-\lambda}(\boldsymbol{D}) \rangle \text{ with } \boldsymbol{A} = \begin{pmatrix} 1+\mu & -\mu & 0 & \cdots & 0 \\ -\mu & 1+2\mu & \ddots & \ddots & \vdots \\ 0 & \ddots & \ddots & \ddots & 0 \\ \vdots & \ddots & \ddots & 1+2\mu & -\mu \\ 0 & \cdots & 0 & -\mu & 1+\mu \end{pmatrix}. \quad (4)$$

ORKA seeks the optimal $\boldsymbol{\lambda}$ by trying to minimize Eq. (4). However, since $\boldsymbol{\lambda} \in \mathbb{N}^N$ this is an integer optimization problem and too hard to solve exactly. Indeed, due to the Lipschitz condition $|\lambda_k - \lambda_{k+1}| \leq C$ there are $(2C+1)^{N-1}$ possible paths $\boldsymbol{\lambda}$ (when fixing $\lambda_1 = 0$). Instead, we solve an approximation of the above term. The exact derivation is quite complicated and can be found in the original work [4]. We only give the brief idea here.

Note that in Eq. (4) the second argument of the inner product is a correlation matrix that compares all columns of $S_{-\lambda}(\boldsymbol{D})$ against each other. All values, no matter how far apart the according measurements were made, influence the final result. ORKA localizes this term, by approximating the first argument $\boldsymbol{A}^{-1}$ by its $K$-bandlimited version

$$\boldsymbol{A}^{-1,[K]} \in \mathbb{R}^{N \times N}, \quad \text{where } \boldsymbol{A}_{jk}^{-1,[K]} = \begin{cases} \boldsymbol{A}_{jk}^{-1} & , |j-k| \leq K \\ 0 & , \text{otherwise} \end{cases}. \quad (5)$$

The coefficients of $\boldsymbol{A}^{-1}$ decay exponentially away from the main diagonal and hence the approximation error is small. ORKA now constructs a graph which we call $K$-approximation graph. Each choice of $\boldsymbol{\lambda}$ corresponds to a path within this graph and the edge weights along this path sum up to $\langle \boldsymbol{A}^{-1,[K]}, (S_{-\lambda}(\boldsymbol{D}))^T S_{-\lambda}(\boldsymbol{D}) \rangle$. Minimizing Eq. (4) can now be approximated by finding the longest path on the $K$-approximation graph which can be done in linear time (for this graph). Because of the locality induced by using the bandlimited inverse (Eq. (5)), the graph size is of order $O((2C+1)^K)$. As $K \ll N$, the complexity of the problem is reduced significantly.

## 2.2 Multiresolution and orthogonal wavelets

Wavelets are specially designed signals that can be used for a localized frequency analysis of a signal $\boldsymbol{v} \in \mathbb{R}^M$. The wavelet transform converts a signal into its lowpass coefficients $\boldsymbol{v}^{\text{low}}$ and highpass coefficients $\boldsymbol{v}^{\text{high}}$. The transform is linear and can be written in matrix form as

$$\boldsymbol{W}\boldsymbol{v} = \begin{pmatrix} \boldsymbol{v}^{\text{low}} \\ \boldsymbol{v}^{\text{high}} \end{pmatrix} \in \mathbb{R}^M. \tag{6}$$

The system matrix $\boldsymbol{W}$ depends on the chosen mother wavelet. For our application, we are especially interested in orthogonal wavelets where $\boldsymbol{W}^{-1} = \boldsymbol{W}^T$, such as Daubechies wavelets [5]. Eq. (6) is often implemented as a filter bank where the highpass coefficients are obtained by a convolution with the wavelet function and the lowpass coefficients by a convolution with the corresponding scaling function. Note that each coefficient vector $\boldsymbol{v}^{\text{low}}, \boldsymbol{v}^{\text{high}} \in \mathbb{R}^{M/2}$ is only half the length of the original signal. This is achieved using a downsampling factor of 2 in the filter bank. While the highpass coefficients contain the details of the signal, the lowpass coefficients can be interpreted as an approximation of the original data with half the resolution. We will use this later in our algorithm to rescale the data into any desired resolution.

Eq. (6) is a one-level wavelet transform. We can apply the same transform again to the low-pass coefficients $\boldsymbol{v}^{\text{low}}$ to obtain the coefficients of the next level. This process can be repeated until only one lowpass coefficient remains. This way, the signal is decomposed into one lowpass coefficient and multiple blocks of highpass coefficients which all correspond to a different level of detail. This decomposition is called multi-resolution analysis of the signal. Although the wavelet transform is not shift-invariant in general, it is according to shifts of even length:

$$\boldsymbol{W}\boldsymbol{v} = \begin{pmatrix} \boldsymbol{v}^{\text{low}} \\ \boldsymbol{v}^{\text{high}} \end{pmatrix} \quad \Rightarrow \quad \boldsymbol{W} S_{2\lambda}(\boldsymbol{v}) = \begin{pmatrix} S_\lambda(\boldsymbol{v}^{\text{low}}) \\ S_\lambda(\boldsymbol{v}^{\text{high}}) \end{pmatrix}. \tag{7}$$

This will be useful in our later analysis.

## 3 multiresolution ORKA

In this section we introduce improved ORKA algorithm. First, we show in what way the resolution of the input data influences the original approach. Afterwards, we show how the multiresolution analysis using wavelets can help overcoming these problems. In the end we give some error bounds of the new approach.

### 3.1 Dependence on data resolution

The resolution of given data always depends on the application and setup of the measurements. Hence, discussing the dependence of ORKA on data resolution is only meaningful with an application in mind. Here, we will demonstrate the dependence for seismic exploration (Fig. 2). For other applications, it follows in a similar manner.

Let $\boldsymbol{D} \in \mathbb{R}^{M \times N}$ be seismic data obtained from $N$ sensors. Each sensor recorded a signal in time containing $M$ data points. The data has a spatial and a time dimension, both having their own resolution. The spatial resolution depends on the positioning of the seismic sensors. For our example, assume the sensors were positioned at locations $x_k$ where $|x_k - x_{k+1}| = 1\text{km}$ for $k = 1, \ldots, N-1$, i.e., we have a spatial resolution of 1km. The time resolution is the time between two samplings. Assume our sensors record 10 samplings per second. Then, the resolution is 0.1s.

The first dependency of ORKA is straightforward. Our object model $S_\lambda(U)$ only allows integer shifts $\lambda \in \mathbb{N}^N$. With the above-mentioned resolution, this is a shift in 0.1s increments. Now assume a seismic wave, modeled by $S_\lambda(U)$, is arriving at one sensor at time $t_k = 2.1$s and at the next sensor at time $t_{k+1} = 2.23$s. The relative difference $|t_k - t_{k+1}| = 0.13$s is not a multiple of 0.1s and can thus not be modeled exactly. Even if the data is accurate enough that the exact time could be extracted, our model is limited by the resolution. We need to increase the resolution to gain more accurate results.

The second dependency of ORKA on the resolution is hidden in the Lipschitz constant $C$. The constant defines how large the relative shift between two neighboring measurements can be. This, however, depends on the application. A seismic wave travels at a speed of at least 5km/s. Hence, for two different sensors at locations $x_k, x_{k+1}$ measuring the same wave at times $t_k, t_{k+1}$ the equation

$$|t_k - t_{k+1}| \leq \frac{|x_k - x_{k+1}|}{5 \text{ km/s}}$$

must hold. With a spatial resolution of 1km, we get $|t_k - t_{k+1}| \leq 0.2$s. For a time resolution of 0.1s, this is equivalent to a shift by 2 samplings, i.e., $|\lambda_k - \lambda_{k+1}| \leq 2$ and $C = 2$. If we increase the time resolution to 0.01s, a difference of 0.2s now corresponds to a shift by 20 samplings and $C = 20$. This increases the runtime of the algorithm a lot since its complexity is $O((2C + 1)^K)$. To keep the runtime low, we need the resolution to be as low as possible. This is in conflict to the high resolution required to obtain accurate shift vectors. Hence, ORKA strongly depends on the resolution of the data.

### 3.2 Multiresolution ORKA

To overcome the resolution dependency of ORKA, we use the wavelet transform to obtain different resolutions of the original data. We start with the lowest resolution to get a starting guess for the shift $\lambda \in \mathbb{N}^N$. Afterwards, the resolution is increased in each step and the shift is updated.

Let $D \in \mathbb{R}^{M \times N}$ be the given data. For $l = 1, \ldots, L$, denote by $D^{\text{low},l} \in \mathbb{R}^{2^{-l}M \times N}$ the lowpass coefficients of the $l$-level wavelet transform. We are only using a column-wise (one-dimensional) wavelet transform here. We use the ORKA algorithm on the lowest resolution $D^{\text{low},L}$ of our data. Our goal is, to keep the complexity of the algorithm low, so we choose the parameter $C = 1$. (We will discuss later, how $L$ has to be chosen such that this choice is correct.) The ORKA algorithm returns our first approximation of the path that we denote $\lambda^L$. In the next step, we want to find $\lambda^{L-1}$, the optimal path for the next higher resolution $D^{\text{low},L-1}$ by using $\lambda^L$. By Eq. (7) we know that the shift vector $\lambda^L$ is equivalent to $2\lambda^L$ for the higher resolution. Hence, we set $\lambda^{L-1} = 2\lambda^L + \lambda'$ where $\lambda'$ is the path update we want to calculate. Note that by Eq. (1) we have $S_{\lambda^{L-1}}(D^{\text{low},L-1}) = S_{\lambda'}(S_{2\lambda^L}(D^{\text{low},L-1}))$ and $S_{2\lambda^L}(D^{\text{low},L-1})$ is already known. We can apply the ORKA algorithm once again, this time with input data $S_{2\lambda^L}(D^{\text{low},L-1})$. We assume that $2\lambda^L$ is already a good approximation of the desired path and the update $\lambda'$ only needs to account for the now higher resolution. The only shift that could not be accomplished in the lower resolution, was a shift by half a sample (one sample in the higher resolution). Thus, we can set $|\lambda'_k - \lambda'_{k+1}| \leq 1$ which results in $C = 1$ again. Now we can repeat this process using $2\lambda^{L-1}$ as approximation for $\lambda^{L-2}$.

This process can be iterated until we obtain $\lambda^0$ which is the shift vector of the original resolution of $D$. If this resolution is still not sufficiently accurate, we can iterate even further. Therefore set $D^{\text{low},0} = D$ and define

$$D^{\text{low},-l} = W^{-1}\begin{pmatrix} D^{\text{low},1-l} \\ 0 \end{pmatrix}, \tag{8}$$

which is the inverse wavelet transform (compare Eq. (6)) where the highpass coefficients are replaced by 0. This is a simple technique to upsample the data and increase the resolution. We can perform another $J \in \mathbb{N}$ iterations to compute $\lambda^{-J}$ as final shift vector. In theory, we can repeat this process indefinitely. However, we are only artificially refining the grid for our shifts. There is no new information added to the data. Depending on the application, the process should saturate after some iterations. Also, note that the shift needs to be rescaled to $\lambda^{-J}/2^J$ in order to fit the original resolution.

Last, we need to find $L$ such that the starting resolution $D^{\text{low},L}$ fits the parameter choice $C = 1$ for the ORKA algorithm. Let $C'$ be the original Lipschitz constant fitting the data $D$ that was obtained from the application setup, e.g., as described in Section 3.1. We require $|\lambda^{-J}_k - \lambda^{-J}_{k+1}| \leq 2^J C'$. The left-hand side can be rewritten as

$$\lambda^{-J}_k - \lambda^{-J}_{k+1} = \sum_{l=-J}^{L} 2^{J+l} \alpha_l \quad \text{with} \quad \alpha_l \in \{-1,0,1\}. \tag{9}$$

The coefficients $\alpha_l$ are derived from the update of each iteration, i.e., $\alpha_l = \lambda'_k - \lambda'_{k+1}$ where $\lambda'$ is the update of the according iteration ($\alpha_L$ is obtained from the starting guess $\lambda^L$). The maximum in Eq. (9) is $2^{J+L+1} - 1$, which is achieved for $\alpha_l = 1$. We obtain

$$2^{J+L+1} - 1 \leq 2^J C' \quad \Leftrightarrow \quad 2^{L+1} - 2^{-J} \leq C' \quad \Leftrightarrow \quad L = \log_2(C' + 2^{-J}) - 1.$$

Since $L \in \mathbb{N}$, we round to result to the next higher number $L = \lceil \log_2(C' + 2^{-J}) \rceil - 1$. Because of this rounding, it is possible that Eq. (9) violates our original Lipschitz condition in the extreme cases. This can be avoided by using a different rescaling that is not a power of 2. However, this requires more complex techniques than the here used wavelet transform. Also note that Eq. (9) is basically a binary representation of the difference that also allows -1 as coefficient. This adds some redundancy to the expression which gives the algorithm more flexibility and makes it more stable against errors. The overall complexity of the new multiresolution ORKA is $O(3^K)$ where we neglect all linear dependencies. Indeed, the runtime depends linear on the data size which means that an excessive oversampling using Eq. (8) will increase the complexity again.

### 3.3 Error bounds

For $\lambda \in \mathbb{N}^N$ choose $\lambda' \in \mathbb{N}^N$, $\lambda'' \in \{\pm 1, 0\}^N$ such that $\lambda = 2\lambda' + \lambda''$. We give a bound on the following term:

$$\left| \langle A^{-1,[K]}, \left(S_{-\lambda}(D)\right)^T S_{-\lambda}(D) \rangle - \langle A^{-1,[K]}, \left(S_{-\lambda'}(D^{\text{low}})\right)^T S_{-\lambda'}(D^{\text{low}}) \rangle \right| \quad (10)$$

$$\leq \left| \langle A^{-1,[K]}, \left(S_{-\lambda}(D)\right)^T S_{-\lambda}(D) \rangle - \langle A^{-1,[K]}, \left(S_{-2\lambda'}(D)\right)^T S_{-2\lambda'}(D) \rangle \right| \quad (11)$$

$$+ \left| \langle A^{-1,[K]}, \left(S_{-2\lambda'}(D)\right)^T S_{-2\lambda'}(D) \rangle - \langle A^{-1,[K]}, \left(S_{-\lambda'}(D^{\text{low}})\right)^T S_{-\lambda'}(D^{\text{low}}) \rangle \right| \quad (12)$$

The first term in Eq. (10) is the equivalent to Eq. (4) using the bandlimited inverse of Eq. (5), i.e., it is the term optimized by the original ORKA method. The second term is the first low-resolution approximation using the lowpass coefficients, i.e., it is the term that will be optimized by a one level ($L = 1$) multiresolution ORKA. The bounds for a higher level ($L > 1$) approach can then be found by applying the bounds given here for each iteration. To get a complete error bound for the multiresolution ORKA approach, we have to add the approximation error of the original ORKA, which was shown to decrease exponentially (see Theorem 7 in [4]).

We first bound Eq. (11). The inner products can be rewritten using Eq. (1) as

$$\left| \langle A^{-1,[K]}, \left(S_{-\lambda}(D)\right)^T S_{-\lambda}(D) - \left(S_{-2\lambda'}(D)\right)^T S_{-2\lambda'}(D) \rangle \right|$$

$$\leq \sum_{j,k=1}^{N} \left| A_{jk}^{-1,[K]} \right| \left| \langle S_{-\lambda_j}(D_{:j}), S_{-\lambda_k}(D_{:k}) \rangle - \langle S_{-2\lambda'_j}(D_{:j}), S_{-2\lambda'_k}(D_{:k}) \rangle \right|$$

$$= \sum_{j,k=1}^{N} \left| A_{jk}^{-1,[K]} \right| \left| \langle D_{:j}, S_{\lambda_j - \lambda_k}(D_{:k}) - S_{2\lambda'_j - 2\lambda'_k}(D_{:k}) \rangle \right|$$

$$= \sum_{j,k=1}^{N} \left| A_{jk}^{-1,[K]} \right| \left| \langle D_{:j}, S_{\lambda_j - \lambda_k}\left(D_{:k} - S_{\lambda''_k - \lambda''_j}(D_{:k})\right) \rangle \right|$$

$$\leq \sum_{j,k=1}^{N} \left| A_{jk}^{-1,[K]} \right| \|D_{:j}\|_2 \left\| D_{:k} - S_{\lambda''_k - \lambda''_j}(D_{:k}) \right\|_2. \quad (13)$$

Since $\lambda''_k - \lambda''_j \in \{\pm 2, \pm 1, 0\}$, this term stays small as long as the data does not change much under small shifts. From Eq. (7) using an orthogonal wavelet transform (i.e., $W^{-1} = W^T$) we can show that

$$\left(S_{-2\lambda'}(D)\right)^T S_{-2\lambda'}(D) = \left(S_{-\lambda'}(D^{\text{low}})\right)^T S_{-\lambda'}(D^{\text{low}}) + \left(S_{-\lambda'}(D^{\text{high}})\right)^T S_{-\lambda'}(D^{\text{high}}).$$

Plugging this into Eq. (12) we obtain the upper bound

$$\left| \langle A^{-1,[K]}, \left(S_{-\lambda'}(D^{\text{high}})\right)^T S_{-\lambda'}(D^{\text{high}}) \rangle \right| \leq \|A^{-1,[K]}\|_F \|D^{\text{high}}\|_F^2. \quad (14)$$

This bound is small as long as the highpass coefficients are small. Eqs. (13) and (14) indicate that, if the data is sufficiently smooth, the error is small. Both bounds hold for all steps of the algorithm (replace $D$ with $D^{\text{low},l}$ and $D^{\text{low}}$ with $D^{\text{low},l+1}$ in Eq. (9)). The right-hand side in Eq. (14) becomes 0 if we upsample the data using Eq. (8).

## 4  Numerical Examples

In all experiments we use Daubechies6 wavelets for the multiresolution approach. For the first experiments we use the data $D \in \mathbb{R}^{512 \times 512}$ with $\alpha, s \in \mathbb{R}$ and

$$D_{jk} = \exp\left(-\left(\frac{mod(j-sk, 512) - 256}{\alpha}\right)^2\right).$$

This is a Gauss kernel where $\alpha$ controls the width and $s$ the shift per column. The modulo operation keeps the kernel inside the data range by making the signal 512-periodic. This data fits our model perfectly whenever $s \in \mathbb{N}$, i.e., we can write it as $D = S_\lambda(U)$ where $\lambda_k = sk$ and $U_{:k} = D_{:1}$. To recover this shift, we need to choose $C = \lceil s \rceil$.

In our first experiment we compare the runtime of ORKA against multiresolution ORKA. Therefore, we choose $\alpha = 10$ and $s = 10$. Since $s \in \mathbb{N}$ both algorithms recover the correct shifts. However, the runtime vastly differs as can be seen in Fig. 1 (left). Here, the runtime for different values of $K$ is plotted in logarithmic scale. Note that the original ORKA approach ran out of memory for $K > 7$. (The space complexity scales in the same way as the runtime.) The experiments were run under MATLAB 2020a on an Intel Core i7-8700 (3.2Ghz) with 64GB memory.

The second experiment uses $\alpha = 10, K = 15$, and a non-integer shift $s = 1/3$. The original ORKA algorithm is not able to reconstruct this shift exactly. We use the multiresolution approach with different levels of upsampling to refine the shift grid. The approximation error $\sum_{k=1}^{512} |\lambda_k - k/3|$ compared to the upsampling factor $2^J$ is shown in Fig. 1 (middle). We also added a line showing the approximation error of the original ORKA algorithm, which returns the optimal integer approximation $\lambda_k = [k/3]$. As we can see, the upsampling increases the accuracy of the shift vector.

Next, we want to test the approximation error of multiresolution ORKA. Therefore, we choose $s = 10$ and $K = 3,15$. Since $s \in \mathbb{N}$, we should be able to recover the shift exactly. We calculate approximation error $\sum_{k=1}^{512} |\lambda_k - sk|$ for different parameters $\alpha$. The results are shown in Fig. 1 (right). For small $\alpha$ the kernel becomes more localized what increases the error bounds (Eqs. (13) and (14)). If the signal is too localized, we can no longer find an accurate shift, independent of the choice of $K$.

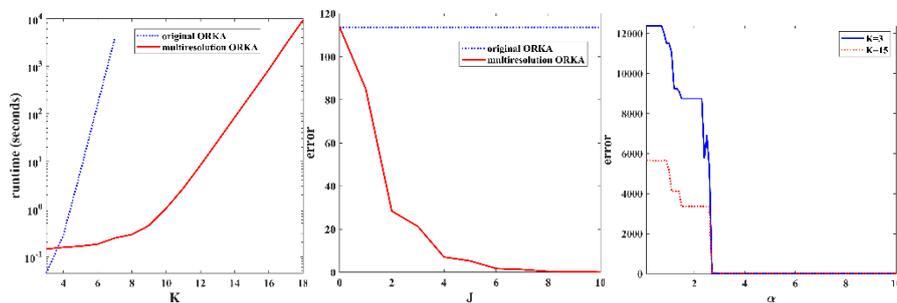

Fig. 1: multiresolution ORKA compared to the original approach: runtime (left), error for non-integer shifts (middle), and error for strongly localized signals (right).

Last, we test our approach on real data from seismic exploration (Fig. 2 (left)). We use $C=3$, $\mu=100$, $K=8$, and an upsampling of $J=3$. Fig. 2 (middle) shows the reconstructed wave movement for the original and the multiresolution ORKA. Because of the upsampling, the new approach is able to catch more details of the movement. This becomes even clearer when zooming into the data (Fig. 2 (right)). Here, the original approach does not detect any shift while the new method finds small movement which can also be seen in the data.

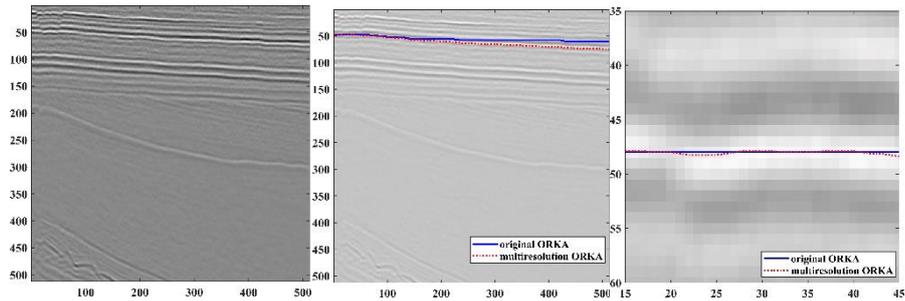

Fig. 2: Seismic exploration data (left) and tracked seismic wave by the original and the multi-resolution approach (middle (full path), right (zoomed in)).

## 5 Conclusion

The original ORKA algorithm can track and extract moving and deforming objects out of given data. However, his runtime and accuracy highly depend on the data resolution. We introduced a multiresolution ORKA approach that uses the wavelet transform to get low resolution approximations of the given data. The object movement is then reconstructed iteratively where in each iteration the resolution is doubled. This way, the new approach can gain the accuracy of the highest resolution used while the runtime is kept low. We have shown that the error compared to the original technique can be bounded as long as the data is sufficiently smooth. Moreover, we demonstrated the advantages of our new approach in numerical experiments.

## References

1. Yurtsever, E., Lambert, J., Carballo, A., Takeda, K. A survey of autonomous driving: Common practices and emerging technologies. *IEEE access*, *8*, 58443-58469 (2020).
2. Rawlinson, N., Hauser, J., Sambridge, M. Seismic ray tracing and wavefront tracking in laterally heterogeneous media. *Advances in geophysics*, *49*, 203-273 (2008).
3. Dutta, K., Krishnan, P., Mathew, M., Jawahar, C. V. Improving CNN-RNN hybrid networks for handwriting recognition. In *2018 16th international conference on frontiers in handwriting recognition (ICFHR)* (pp. 80-85). IEEE (2018).
4. Bossmann, F., Ma, J. ORKA: Object reconstruction using a K-approximation graph. *Inverse Problems, accepted (available online* https://doi.org/10.1088/1361-6420/aca046*)*, (2022).
5. Daubechies, I. *Ten lectures on wavelets*. Society for industrial and applied mathematics (1992).